\documentclass[a4paper,11pt]{article}
\usepackage{latexsym,amssymb}
\newtheorem{prop}{Proposition}[section]
\newtheorem{thm}[prop]{Theorem}
\newtheorem{Def}[prop]{Definition}
\newtheorem{lemma}[prop]{Lemma}

\newcommand{\F}{\mathbb F}
\newcommand{\Cr}{\mathrm {Cr}}
\newcommand{\s}{\mathcal  L}
\newcommand{\gl}[2]{\mathrm{GL}(#1,#2)}
\newcommand{\pf}{\noindent{\bf Proof~}}

\begin{document}
\title{\huge{Hirsch-Plotkin radical of  stability groups}}
\author{Carlo Casolo  \\ \small
     Dipartimento di Matematica ``U.\,Dini",
Universit\`a di Firenze\\ \small
    Viale Morgagni 67A, \ I-50134 Firenze, Italy.\\[-.5ex] \small
     \texttt{casolo\,@\,math.unifi.it}\\[1ex]
\and Orazio Puglisi\\ \small
Dipartimento di Matematica ``U.\,Dini",
Universit\`a di Firenze\\ \small
Viale Morgagni 67A, \ I-50134 Firenze, Italy\\ \small
\texttt{puglisi\,@\,math.unifi.it}}
\date{}
\maketitle

%
%
%%%%%%%%%%%%%%%%%%%%%%%%%%%%%%%%%%%%%%%%%%%%%
%
%
\section{Introduction}
\bigskip
Let $V$ be a vector space over the field  $\F$, and  $G$ a subgroup of $\gl{V}{\F}$. The group  $G$ is said to be \emph{unipotent} if , for each $g\in G$, there exists $n=n(g)\in \mathbb N$, such that $(g-1)^n=0$ in $\mathrm{End}(V)$. When $V$ has finite dimension, the unipotent subgroups of $\gl{V}{\F}$ are well known. It is easy to prove in that case that a subgroup  $G$ of $\gl{V}{\F}$ is unipotent if and only if it stabilizes a series in $V$. Hence unipotent groups are subgroups of \emph{series stabilizers}. From this fact it readily follows that unipotent subgroups are nilpotent of class at most $\dim(V)-1$, and they are torsion-free or $p$-groups according to characteristic of $\F$ being $0$ or a prime $p$. Moreover several other insights about the structure of unipotent groups can be drawn using their action on the natural module. When the vector space $V$ has infinite dimension, the concept of series must be considered with care. Although it is easy to figure what an 
 \emph{ascending} or \emph{descending} series should be, the definition of \emph{series} can be given in different ways. We shall give the precise definition (at least the one we find more suitable for us) in the next section, but remark that, once this has been done, it is straightforward to define stability groups also in this more genereal setting. However, when the series is infinite,  unipotent groups and series stabilizers do not  coincide anymore. It is easy to produce examples of series stabilizers containing non-unipotent elements and, on the other hand, there exist groups acting unipotently which are  not contained in any stability group. We shall remind examples of these situations in the following section. 

Although we do not expect that the elements of a stability group $S(\mathcal L)$ act unipotently when the series $\mathcal L$ is infinite, it is conceivable that this happens if we consider some particular subgroup. A natural candidate is the Hirsch-Plotkin radical of $S(\s)$, and we prove that indeed its elements are always unipotent in their action on the natural module. 
This fact suggests that the structure of the Hirsch-Plotkin radical of a stability group may be succesfully investigated. However this hope is a little bit too optimistic, because some problems arise when the series $\mathcal L$ contains  certain kinds of descending chains.  Nevertheless the case of vector spaces of countable dimension can be fully described.

\

{\bf Theorem A}\emph{ Let $V$ be an $\F$-vector space of dimension at most $\aleph_0$, $\mathcal L$ a series in $V$ and $H(\s)$ the Hirsch-Plotkin radical of $S(\mathcal L)$, the stability group of $\mathcal L$. Then $H(\s)=\{ g\in S(\mathcal L) \mid g$ stabilizes a finite subseries of $\mathcal L\}$. Moreover $H(\s)$ is a Fitting group.
}

\

We are then lead to state the following conjecture

\

{\bf Conjecture }\emph{ Let $V$ be an $\F$-vector space, $\mathcal L$ a series in $V$ and $H(\s)$ the Hirsch-Plotkin radical of $S(\mathcal L)$, the stability group of $\mathcal L$. Then $H(\s)=\{ g\in S(\mathcal L) \mid g$ stabilizes a finite subseries on $\mathcal L\}$.}

\

We are able to prove that the conjecture holds in several situations, but the general case  remains unsolved.

If the conjecture holds for $S(\s)$, the Hirsch-Plotkin radical has the following strong property:
\begin{itemize}
\item for every $g\in H(\s)$ the group $g^{S(\s)}$ is nilpotent.
\end{itemize}

This is easily seen because $g^{H(\s)}\leq H(\s)$ stabilizes the same subseries as $g$.

\

A related result can be deduced using work of Roseblade. The results proved in \cite{blade} imply that, when the vector space $V$ has an $\s$-adapted base (see  section 3 for the definition), then the  Baer radical of $S(\s)$ coincides with the Fitting subgroup.

\

 At this stage it is natural to ask whether every Fitting group can be embedded into  the Hirsch-Plotkin radical of a suitable stability group or, at least, if it can be faithully represented as a unipotent subgroup of a stability group.
 We consider this problem in the last section of this paper. It can be easily seen that this question makes sense only for groups which are  torsion-free  or $p$-groups. We show that such groups admit unipotent representations satysfing a rather strong condition. The existence of this representation is then used to prove the following theorem.
   
 \
 
 {\bf Theorem B}\emph{ Let $G$ be a  Fitting group which is torsion-free or a $p$-group. Assume that there exists a countable subset $X\subseteq G$ such that $G=X^G$. Then there exist an $\F$-vector space $V$ and a series $\s$ in $V$, such that $G$ can be embedded into the Hirsch-Plotkin radical of $S(\s)$.
}

\

The above theorem applies, e.g., to countable groups, but we have been unable to extend it to groups of arbitrary cardinality.  

 It is worth remarking that, for every cardinal $\kappa$, there exists a Fitting group $G$ of cardinality $\kappa$, such that $G=X^G$ for a suitable countable subset $X$ of $G$. Moreover $G$ can be chosen to be a $p$-group or torsion-free. An example of such group will be described in the last section.

 It might be of interest to point out that the problem of embedding groups, satisfying some nilpotency condition, into series stabilizers has been investigated by  several authors. In \cite{w2} Wehrfritz shows that every nilpotent group in embeddable in a group of unitriangular matrices, defined over a suitable divison ring, while Leinen discusses the representability of Fitting $p$-groups as unipotent groups of finitary transformations (see \cite{leinen}). 
 \bigskip

\section{ Notation and definitions}

In infinite dimensional spaces, the notion of series must be considered with some care. The following definition is the one best suited for us.
\

\begin{Def}\label{series} Let $\F$ be a field and $V$ an $\F$-vector 
space. A set $\s$ of subspaces of $V$ is said to be  a \emph{ series} in $V$ if
\begin{enumerate}
\item $0$ and $  V$ belong to $\s$;
\item the set $\s$ is linearly ordered;
\item for every $\mathcal F\subseteq \s$, both $\cap\{ W\mid W\in \mathcal 
F\}$ and $\cup\{ W\mid W\in \mathcal 
F\}$ belong to $\s$.
\end{enumerate}
A \emph{subseries} of $\s$ is any subset of $\s$ which is still a series. 
In particular every finite subset of $\s$ containing $0$ and $V$ is a 
subseries.
\end{Def}

\

When the series is finite, we define its \emph{length} as the number of its non-trivial elements.

\

In many cases we shall be looking at quotients of the form $W/U$ for some $U,W\in \s$. We call such a quotient a \emph{section} of $\s$.

If $\s=\{V_\lambda\mid \lambda\in \Lambda\}$ is a series in $V$, then the 
set $\Lambda$ can be endowed by a natural order by setting $\lambda\leq 
\mu$ if and only if $V_\lambda\leq V_\mu$. The set $(\Lambda, \leq)$ is 
linearly ordered, and each non empty subset of $\Lambda$ has both an 
infimum and a supremum.
 
\

\begin{Def} Let $\s$ be a series in the $\F$-vector space  $V$. A 
\emph{jump} of  $\s$ is an ordered pair $(B,T)$ of elements of $\s$ such that:
\begin{itemize}
\item $B<T$;
\item if $U\in \s$ and  $B\leq U\leq T$, then $U\in \{B, 
T\}$.
\end{itemize}
\end{Def}

Jumps are easy to produce. Given a series $\s$ choose a non-trivial vector 
$v\in V$ and let $B_v=\cup\{ U\in \s\mid v\not\in U\}$, 
$T_v=\cap\{ U\in \s\mid v\in U\}$. The pair $(B_v,T_v)$ is a jump for $\s$ 
and any jump can be obtained in this way.  We say that $v$ \emph{belongs }
to the jump $(B,T)$,  if $(B,T)= (B_v,T_v)$. This happens if and only if $v\in T\setminus B$ and it is clear that a non-trivial 
vector $v$ belongs to just one jump of $\s$.

\

Let $\mathcal B$ be a basis of $V$ and, for each jump $j=(B,T)$ define $\mathcal B_j=\{v\in \mathcal B \mid v \textrm{ belongs to } j\}$. The basis $\mathcal B$ is said to be \emph{$\s$-adapted} if, for  each jump $j=(B,T)$, the set $\mathcal B_j+B=\{v+B\mid v\in \mathcal B_j\}$ is a basis of $T/B$. 
It is worth remarking that the existence of $\s$-adapted basis is not granted in infinite dimension. Perhaps the easiest example is the following.

\

Choose any field $\F$ and let $V=\F^\omega$. For each $i\in \omega $ set $V_i=\{(v_j)\mid v_j=0 \, \, \forall j<i\}$. The series $\s=\{V_i\mid i\in \omega\}\cup \{0\}$ has countably many jumps, each of dimension one, while $V$ has uncountable dimension. Hence no $\s$-adapted basis can  exists.

\

Given any $W\in \s$, two series can be  defined in $W$ and $V/W$. 
We set $\s\cap W=\{ U\in 
\s\mid U\leq W\}$ and $\s/W=\{ U/W\mid U\in \s, \, W\leq U\}$.

\

The set of series in a fixed vector space,  endowed with its natural order,  is an inductive set, whence   every series $\s$ can be extended to a maximal one, namely a series $\overline{\s}$ whose jumps have  dimension $1$.

\begin{Def}Let $V$ be an $\F$-vector space, consider a series $\s$ and its 
set of jumps $\mathcal J$. We define
$$
S(\s)=\{ g\in \gl{\F}{V}\mid [T,g]\leq B \quad \forall \, (B,T)\in 
\mathcal J\}
$$
\end{Def}
The set $S(\s)$ is  a subgroup of $\gl{V}{\F}$ called  the 
\emph{stability group of } $\s$ or the \emph{stabilizer} of $\s$.

If $V$ has finite dimension $n$, then any series $\s$ in $V$ is finite, and 
$S(\s)$ is isomorphic to a subgroup of $\mathrm{Tr}(n,\F)$, the group of $n\times n$ 
unitiangular matrices over $\F$. Therefore $S(\s)$ is nilpotent, and its 
nilpotency class is bounded by $n-1$. However we should not expect that 
$S(\s)$ satisfies any kind of solubility condition when the series $\s$ 
contains infinitely many elements, as the following example shows.

\

{\bf Example 1} {\it  Every  countable free group  faithfully stabilizes an ascending 
series in a 
suitable vector space.}

\

Let $G$ be any  free group, choose a prime $p$ and let 
$\F_p$ denote the field with $p$ elements. It is 
well known that $G$ is residually a finite $p$-group, and we can find a 
descending chain $\{N_i\mid i\in \omega\}$ of normal subgroups such that 
$G/N_i$ is a $p$-group for every $i$, and $\cap_{i\in \omega}N_i=1$. For 
each $i$ let $V_i$ be a faithful $\F_p(G/N_i)$-module of finite dimension and set 

$$V=\mathrm{Dir}\{V_i\mid i\in \omega\}.$$ 
The group $G$ acts faithfully on 
$V$ and normalizes the subspaces 
$$W_n=\mathrm{Dir}\{V_i\mid n\geq i\in \omega\}$$
Since $W_{n+1}/W_{n}$ is a finite $\F_p$-vector space and   
$G/C_G(W_{n+1}/W_{n})$ 
is a 
finite $p$-group, $G$ stabilizes a finite series in each factor 
$W_{n+1}/W_{n}$. The preimages of the members of these series, together with 
$0$, form an ascending  series $\s$  in $V$ which is stabilized by $G$. Thus $S(\s)$ 
contains a copy of $G$. The reader can easily  find in $V$ a \emph{descending} series stabilized by $G$. 

\

Another well known property of stability groups in finite dimension, is 
their unipotency. We recall that an element $g\in \gl{\F}{V}$ is said to 
be unipotent if there exists $n\in \mathbb N$ such that $(g-1)^n=0$ in the 
ring $\mathrm{End}_\F(V)$. When $V$ has finite dimension, a subgroup 
$G\leq \gl{\F}{V}$ stabilizes a series in $V$  if and only if each element 
of $G$ is unipotent. Actually this happens if and only if there is an 
$n\in \mathbb N$ such that $(g-1)^n=0$ for every element of $G$.  
On the other hand  stability groups of infinite series may not be unipotent.
E.g. consider a vector space $V$ of countable dimension over any field 
$\F$ and select a basis $B=\{ v_i\mid i\in \mathbb Z\}$. There is an 
element $g\in \gl{V}{\F}$ such that $v_ig=v_i+v_{i-1}$. It is clear that 
$g$ stabilizes a series in $V$, but $(g-1)^n\not=0$ for all $n\in \mathbb 
N$.

Moreover groups acting unipotently may not stabilize any series. 

\

{\bf Example 2} {\it Unipotent actions of Tarski groups.}
 
 \

 Let $G$ be a Tarski group. Then $G$ is infinite and there exists a prime $p$ such that every proper non-trivial subgroup of $G$ is cyclic of order $p$. Let $\F$ be any field of characteristic $p$ and set $V=\F G$ for the group algebra of $G$ over $\F$. The group $G$ acts naturally on $V$. Moreover, for each $g\in G$, we have $(g-1)^p=g^p-1=0$ in $\mathrm{End}(V)$, so that the action of $G$ on $V$ is faithful and unipotent. Should $G$ stabilize a series in $V$ then, by Lemma 4 of \cite{HH}, $G$ would have a series with abelian factors, and this is impossible.

 \

\section{ Hirsch-Plotkin radical of  a stability group}
\bigskip

This section is devoted to the investigation of the Hirsch-Plotkin radical of a stability group. The first result we prove concerns the unipotency of its elements. Some lemmata are needed.

Let $U$ be any element in $\mathcal L$. For every $\phi\in \mathrm{Hom}(V/U,U)$ define $x_\phi: V\longrightarrow V$ by $vx_\phi= v+(v+U)\phi$. It is readily seen that $x_\phi$ belongs to $G=S(\s)$ and actually $A_U=S(\{0,U,V\})=\{ x_\phi\mid \phi\in \mathrm{Hom}(V/U, U)\}$. The group $A_U$ is abelian and normal in $G$, and the action of $G$ on $A_U$ can be described very clearly. In fact a straightforward calculation shows that the map
\begin{displaymath}
\begin{array}{rr}
\mathrm{Hom}(V/U, U) & \longrightarrow A_U\\
\phi & \longmapsto x_\phi
\end{array}
\end{displaymath}
is an isomorphism of $G$-modules.

\

\begin{lemma}\label{lemma 1}Let $V$ be an infinite dimensional $\F$-vector space and $g\in \mathrm{ GL}(V,\F)$ a unipotent element. Then the dimension of $V/[V,g]$ is infinite
\end{lemma}

\pf The element $g$ is unipotent, so that there exists $k\in \mathbb N$ such that $(g-1)^k=0$. Call the \emph{exponent} of $g$ the minimal $n$ such that $(g-1)^n=0$
If the exponent of $g$ is $1$ then $g$ is the identity and  the claim is trivially true. Assume the lemma holds for unipotent elements of exponent smaller that $n$, and let $g$ have exponent $n$. Set $W=V/[V,_{n-1} g]$ and consider the transformation $x$ induced by $g$ on $W$. If $W$ has infinite dimension we use inductive hypothesis to see that $W/[W,x]$ has infinite dimension. From this the claim follows. Hence $W$ has finite dimension. But then $[V,_{n-1} g]$ has finite codimension and, a fortiori, this holds for $[V,g]$. The subspace  $[V,g]$ is isomorphic to $V/\ker(g-1)$, whence it has finite dimension, because $[V,_{n-1} g]\leq \ker(g-1)$. On the other hand $[V,_{n-1} g]\leq [V,g]$ hence, being $[V,_{n-1} g]$ of finite dimension and codimension, $V$ is finite dimensional. This contradiction proves the claim.\hfill $\Box$

\

\begin{lemma}\label{lemma 2}Let $V$ be an $\F$-vector space, $\s$ a series in $V$ and $U$ any element of $\s$. Choose $\phi\in \mathrm{Hom}(V/U,U)$, $t\in G=S(\s)$ and assume that $[V,t]+U/U\leq \ker(\phi)$. Then, for every $v\in V$ and $1\leq k\in \mathbb N$,  $v[x_\phi,_k t]=v+(v+U)\phi (t-1)^k$.
\end{lemma}

\pf We argue by induction on $k$, since the case $k=1$ is clearly true. For every $v\in V$  write 
$$v[x_\phi, _k t]=v[x_\phi,_{k-1} t]^{-1}[x_\phi,_{k-1} t]^t
=(v-(v+U)\phi (t-1)^{k-1})[x_\phi,_{k-1} t]^t$$
Since $(v+U)\phi (t-1)^{k-1}$ is in $U$, it is left fixed by the transformation 
$[x_\phi,_{k-1} t]^t$. The  inductive hypothesis yelds
$$
v[x_\phi, _k t]=v[x_\phi,_{k-1} t]^t-(v+U)\phi (t-1)^{k-1}=
$$
$$
(vt^{-1}+(vt^{-1}+U)\phi(t-1)^{k-1})t-(v+U)\phi (t-1)^{k-1}=
$$
$$
v+(vt^{-1}+U)\phi(t-1)^{k-1}t-(v+U)\phi (t-1)^{k-1}
$$
Since $vt^{-1}=v+vt^{-1}(1-t)$ and $[V,t]+U/U\leq \ker(\phi)$, we get
$$
v+(vt^{-1}+U)\phi(t-1)^{k-1}t-(v+U)\phi (t-1)^{k-1}=
$$
$$
v+(v+U)\phi(t-1)^{k-1}t-(v+U)\phi (t-1)^{k-1}=
v+(v+U)\phi(t-1)^k
$$
and the claim is proved.\hfill $\Box$

 \
 
 Let us consider a series $\s$ in the vector space $V$. This series can be extended to a maximal one $\overline{\s}$, and it is readily seen that $S(\s)$ is normal in $S(\overline{\s})$. Thus the Hirsch-Plotkin radical of $S(\s)$ is contained in the Hirsch-Plotkin radical of $S(\overline{\s})$ and, for this reason, there is no loss of generality, for our purposes, in assuming that the series $\s$ is alway maximal.
 
 \

\begin{thm}\label{1}Let $V$ be an $\F$-vector space, $\mathcal L$ a series in $V$ and $g\in G=S(\mathcal L)$ a non-unipotent element.Then $g$ does not belong to the the Hirsch-Plotkin radical of $G$.
\end{thm}

\pf By 12.3.2 of \cite{robinson}, the Hirsch-Plotkin radical of $G$ consists of left-engel elements. In order to prove our claim we show that a non unipotent element $g\in G$ can not be left-engel. Namely we show that, for every $n\geq 1$ there exists $x=x(n,g)$ such that $[x,_ng]\not=1$.  Several cases should be considered.

\

{\bf Case 1.}\emph{ There exists $0\not=U\in \mathcal L$ such that $g$ is not unipotent in its action on $V/U$.}

 Choose a jump $(B,T)$ with $T\leq U$ and consider the series $\mathcal S=\s/B$. The stabilizer of this series is an homomorphic image of $G$, so that the claim will be proved if we show that the transformation induced by $g$ on $V/B$ does not belong to the Hirsch-Plotkin radical of the stabilizer of $\mathcal S$. Without loss of generality we can therefore assume that $\mathcal L$ has a first non-trivial element $U=\langle u\rangle$ and $g$ is not unipotent in its action on $V/U$. It is then possible to find a $\phi\in \mathrm{Hom}(V/U, U)$ such that $[V,_n g]+U/U$ is not contained in $\ker(\phi)$ for all $n\in \mathbb N$. It is easily seen that $z_n=[x_\phi,_n g]$ acts as $vz_n=v+(v(g^{-1}-1)^n+U)\phi$. Thus $1\not=z_n$ for all $n\in \mathbb N$. 
 
 \

{\bf Case 2.}\emph{ For each $0\not=U\in \s$, the element $g$ acts unipotently on $V/U$}

\

The analysis of this situation is divided into some subcases.

\

{\bf Case 2.1.} \emph{ There exists $0\not=U\in \s$ of infinite codimension}

Lemma \ref{lemma 1} can be applied to the action of $g$ on $V/U$ to show that  $V/[V,g]+U$ has  infinite dimension. For each $n\in \mathbb N$ select $u_n\in U$ in such a way that $[u_n, _ng]\not=0$. This is indeed possible since $g$ is not unipotent on $V$ while it acts unipotently on $V/U$. Choose now $\phi\in \mathrm{Hom}(V/U, U)$ in such a way that $[V,g]+U\leq \ker(\phi) $ and $\langle u_n\mid n\in \mathbb N\rangle\leq \mathrm{Im}(\phi)$. Such an homomorphism always exists because $V/[V,g]+U$ is infinite dimensional. Lemma \ref{lemma 2} can now be invoked to see that $z_n=[x_\phi, _n g]$ is not $1$. 

\

{\bf Case 2.2.} \emph{ The space $V/U$ has finite dimension for all $0\not=U\in \s$ but there exists $w\in V$ such that $w(g-1)^n\not=0$ for all $n\in \mathbb N$.}

 The series $\s$ has a maximal proper member, say $W$, and we may assume, without loss of generality, that $w\in W$. Choose $\phi\in \mathrm{Hom}(V/W, W)$ and consider $x_\phi$. An application of lemma \ref{lemma 2} shows that $z_n=[x_\phi,_ng]$ acts as $vz_n=v+(v+W)\phi(g-1)^n$. If $(u+W)\phi=w$ we deduce that, for every $n$,  $uz_n\not=u$, so that $z_n\not=1$.
 
 \

{\bf Case 2.3.} \emph{ The space $V/U$ has finite dimension for all $0\not=U\in \s$ and, for each $v\in V$, there exists $k=k(v)$ with $v(g-1)^k=0$.}

\

 Without loss of generality we assume that $\s$ is maximal. It is readily seen that its order
 type is the reverse order on $\omega+1$, hence $\s=\{V_i\mid i\in \omega\}\cup \{0\}$ and 
 $\dim(V/V_i)=i$ for all $i\in \omega$.

Set $K_n=\ker(g-1)^n$. Our hypothesis imply that $V=\bigcup_{n\in \mathbb N}K_n$. 
In order to make the proof more transparent, we single out some properties of the chain $\{K_i\mid i\in \omega\}$.
\begin{itemize}
\item[(1)] For every $n\in \omega$ there exists $i\in \omega$ such that $K_n+V_i<K_{n+1}+V_i$.
\end{itemize}

 Should this not be true, we would have 
 $K_n+V_i=K_{n+1}+V_i$ for all $i\in \omega$, so that $K_{n+1}=K_{n+1}\cap (K_n+V_i)=K_n+(K_{n+1}\cap V_i)$. 
 Thus $$K_{n+1}(g-1)^n=[K_n+(K_{n+1}\cap V_i)](g-1)^n=(K_{n+1}\cap V_i)(g-1)^n\leq V_{i+1}.$$
  Since this holds for all $i$, we have $K_{n+1}(g-1)^n\leq \cap_{i\in \mathbb N} V_i=0$ showing that $K_{n+1}\leq K_n$, a contradiction.
  
  \begin{itemize}
\item[(2)] For every $i,n\in \omega$ $V_i$ is not contained in $K_n$.
\end{itemize}

 If $K_n\cap V_i=V_i$ then $V_i\leq K_n$ and $K_n$ turns out to have finite codimension $d$. It is then an easy matter to show that $(g-1)^{n+d}=0$, a contradiction. Hence $V_i$ contains elements not lying in $K_n$. 
 
 \begin{itemize}
\item[(3)] If $V_i+K_n<V$, then $V_i+K_n+[V,g]<V$.
\end{itemize}
 
 The space $V_i+K_n$ is $g$-invariant, so that $g$ acts on the finite dimensional space $W=V/(V_i+K_n)$ and, being its action unipotent, we have $[W,g]<W$. Hence $V/(V_i+K_n+[V,g])$ is non trivial. 
 
 \

 \begin{itemize}
 \item[(4)] Fix   $i,n\in \omega$. There exists $n_0$ such that   $K_m\cap V_i>K_n\cap V_i$ for all $m\geq n_0$.
  \end{itemize}
  
  If $K_m\cap V_i=K_n\cap V_i$ for all $m\geq n$, we have 
 $$
 V_i\cap K_n= \bigcup_{n\leq m}V_i\cap K_m\leq V_i\cap(\bigcup_{n\leq m}K_m)=V_i\cap V=V_i
 $$
 By point $(2)$ this can not hold, so $K_m\cap V_i>K_n\cap V_i$ for all sufficiently large $m$.

 \begin{itemize}
 \item[(5)] Fix   $i,n\in \omega$. If $V_i+K_n+[V,g]<V$, there exists $n_0$ such that   $V_i+K_n+[V,g]<V_i+K_m+[V,g]$ for all $m\geq n_0$.
  \end{itemize}
 
 If $V_i+K_n+[V,g]=V_i+K_m+[V,g]$ for all $m\geq n$, then 
 $$
 V_i+K_n+[V,g]=\bigcup_{m\geq n}V_i+K_m+[V,g]\supseteq \bigcup_{m\geq n}K_m=V
 $$
 
 From $(4)$ and $(5)$ it follows that 
 
 \begin{itemize}
 \item[(6)] Given   $i,n\in \omega$ such that  $V_i+K_n+[V,g]<V$, there exists $n_0$ such that   $K_m\cap V_i>K_n\cap V_i$ and $V_i+K_n+[V,g]<V_i+K_m+[V,g]$ for all $m\geq n_0$.
  \end{itemize}

We construct an element $x$ in $G$ such that $\langle g^x, g\rangle$ is not nilpotent.

Start by setting $A_1=V_1$ and choose $i_1$ in such a way that $A_1+K_{i_1}<V$. Set $B_1=K_{i_1}$.  By $(6)$ it is possible to find $B_2\in \{ K_n\mid n\in \omega\}$ satisfying $B_2\cap A_1>B_1\cap A_1$ and $B_2+A_1+[V,g]>B_1+A_1+[V,g]$. 
Find    $\psi_1\in \mathrm{Hom}(V/(A_1+B_1+[V,g]),A_1)$ with
\begin{itemize}

\item $\mathrm{Im}(\psi_1)\leq B_2\cap A_1$ but $\mathrm{Im}(\psi_1)$ is not contained in  $B_1\cap A_1$;
\item $(B_2+A_1+[V,g]/B_1+A_1+[V,g])\psi_1$ is not contained in $B_1\cap A_1$.
\end{itemize}   

Say we have already selected $\{A_i\mid i=1, \dots , m-1\}\subseteq \s$, $\{B_i\mid i=1, \dots , m\}\subseteq \{K_n\mid n\in \omega\}$, and found, for each $i=1, \dots , m-1$ an homomorphism $\psi_i\in \mathrm{Hom}(V/(A_i+B_i+[V,g]),A_i)$, satisfying
\begin{itemize}
\item[(a)] $A_i>A_{i+1}\not=0$ for all $i=1, \dots , m-2$;
\item[(b)] $B_{i+1}\cap A_i>B_i\cap A_i$ and $(A_i+B_{i+1}+[V,g])>(A_i+B_i+[V,g])$, for all $i=1, \dots , m$;
\item[(c)] for all $i=1, \dots , m$,  $\mathrm{Im}(\psi_i)\leq B_{i+1}\cap A_i$ but $\mathrm{Im}(\psi_i)$ is not contained in $B_i\cap A_i$;
\item[(d)] $(B_{i+1}+A_i+[V,g]/B_i+A_i+[V,g])\psi_i$ is not contained in $B_i\cap A_i$.
\end{itemize}

\

To enlarge this set of data define $A_m$ to be any non-trivial element of $\s$ for which $A_m<A_{m-1}$ and $A_m+B_m<V$. The existence of $A_m$ is ensured by $(1)$. Use $(6)$ to find $B_{m+1}\in \{K_n\mid n\in \omega\}$ with $B_{m+1}\cap A_m>B_m\cap A_m$ and 
$B_{m+1}+A_m+[V,g]>B_m+A_m+[V,g]$. The existence of $\psi_m\in \mathrm{Hom}(V/(A_m+B_m+[V,g]),A_m)$, fulfilling conditions $(c)$ and $(d)$ is clear.

\

By iterating this procedure we end up with the following data: 
 two  sets of subspaces 
 $\{A_i\mid i\in \omega\}\subseteq \{V_i\mid i\in \omega\}$, 
 $\{ B_i\mid  i\in \omega\}\subseteq \{ K_i\mid  i\in \omega\}$ and a set of homomorpsisms $\{\psi_i\mid \psi_i\in \mathrm{Hom}(V/(A_i+B_i+[V,g]),A_i)\}$, satisfying the  conditions from $(a)$ to $(d)$ for all $i\in \omega$.

\

For each $i$ let $\pi_i:V\longrightarrow V/(A_i+B_i+[V,g])$ be the canonical projection and define $\phi_i=\pi_i\psi_i$. Given any $v\in V$, the set $\{ B_i\mid  v\not\in B_i\}$ is finite, so that we can define the endomorphism $\eta=\sum_{i\in \mathbb N} \phi_i$, because each vector $v$ is  in $\ker(\phi_i)$ almost always. We show that $\eta^2=0$. For each $i$ we have $\mathrm{Im}(\phi_i)\leq A_i\cap B_{i+1}\leq B_{i+1}\leq \ker(\phi_{i+1})$. Since ${\phi_i}^2=0 $ we get  
$\eta^2=\sum_{i>j}\phi_i\psi_j$. On the other hand, when $i>j$, $A_i\leq A_j$, hence 
$\mathrm{Im}(\phi_i)\leq B_{i+1}\cap A_i\leq A_i\leq A_j\leq \ker(\phi_j)$, so that $\phi_i\phi_j=0$ and $\eta^2=0$ as claimed. Thus the element $x=1+\eta$ belongs to $\mathrm{GL}(V, \mathbb F)$ and it is readily seen to lie in $G$. Since $g\phi_i=\psi_i$, it is  easy to  show that $[x,g]=1+\eta(g-1)$. We  prove that the relation $[x,_n g]=1+\eta(g-1)^n$ holds for every $n\geq 1$. Since the claim is true for $n=1$, we use induction on $n$ and assume this fact holds for $n-1$. Then $[x,_n g]=[[x,_{n-1} g],g]=[1+\eta(g-1)^{n-1},g]$. The square of the  endomorphism  
$\eta(g-1)^{n-1}$ is $0$ and we get
$$
[1+\eta(g-1)^{n-1},g]=(1-\eta(g-1)^{n-1})(1+\eta(g-1)^{n-1})^g=
$$
$$
(1-\eta(g-1)^{n-1})(1+\eta(g-1)^{n-1}g)=1+\eta(g-1)^n
$$

To show that $\langle g^x,g\rangle$ is not nilpotent, it is sufficient to prove that none of  the commutators $[x,_n g]$ is trivial. Choose $n\geq 1$ and pick $i$ in such a way that $K_n\leq B_i$. If $B_i=K_r$ we prove that $[x,_n g]\not=1$ by showing that $[x,_r g]$ is non trivial.
By the above calculation we have 
$$[x,_r g]=1+\eta(g-1)^r=1+\sum_{j\geq i}\phi_j(g-1)^r
$$
because $\mathrm{Im}(\phi_l)\leq \ker((g-1)^r)=B_i$ for all $l\leq i$. Let $v\in B_{i+1}$ be such that $(v)\phi_i\in B_{i+1}\setminus B_i$. Since $(v)\phi_j=0$ for all $j>i$, we get 
$$
(v)[x,_r g]=v+\sum_{j\geq i}(v)\phi_j(g-1)^r=v+(v)\phi_i(g-1)^r
$$
But $(v)\phi_i$ does not belong to $B_i=K_r=\ker((g-1)^r)$, whence $(v)[x,_r g]\not=v$, showing that $[x,_r g]$ is not $1$ and proving the claim. \hfill$\Box$ 

\bigskip 

Theorem  \ref{1} can be stated as follows
\

\begin{thm}\label{first theorem} Let $V$ be an $\F$-vector space, $\s$ a series in $V$ and $S(\s)$ its stabilizer. Then the elements of the Hirsch-Plotkin radical of $S(\s)$ are unipotent.
\end{thm}

\bigskip

Given any series $\s$ in the vector space $V$, we can consider the set $F(\s)=\{g\in S(\s)\mid g$ stabilizes a finite subseries of $ \s\}$. It is immediate to show that $F(\s)$ is a normal subgroup of $S(\s)$ and that, for every $g\in F(\s)$, the group $g^{S(\s)}$ stabilizes the same subseries stabilized by $g$. In particular  $g^{S(\s)}$ is nilpotent whenever $g$ belongs to $F(\s)$  and   $F(\s)$ is contained in the Fitting radical of $S(\s)$. A fortiori $F(\s)$ is a subgroup of  $H(\s)$. We conjecture that $H(\s)=F(\s)$ but we have been unable to prove that this equality holds in full generality. It is  however possible to show that, in some relevant cases, our conjecture holds.  The first situation we discuss is the case of vector spaces of  countable dimension.
The proof relies heavily on the following lemma.

\bigskip

\begin{lemma}\label{decomposition} Let $V$ be an $\F$-vector space of countable dimension, and $\s=\{ V_i\mid i\in \omega\}\cup\{0\}$ a descending series of subspaces. Then there exist subspaces $\{A_i\mid i\in \omega\}$ and a strictly increasing map $\sigma: \omega \longrightarrow \omega$, such that
\begin{enumerate}
\item $V=\bigoplus_{i\geq 1}A_i$;
\item $A_i\simeq V_{\sigma(i-1)}/V_{\sigma(i)}$ for all $i\geq 1$.
\end{enumerate}  
\end{lemma}
 
 \pf Fix a basis $\mathcal B=\{ v_i\mid i\in \omega\}$ and, for each $i\geq 1$, let $C_i$ be any complement to $V_i$ in $V$ chosen in such a way that $C_i\leq C_{i+1}$ for all $i\in \omega$. Set $L_i=\langle C_i, v_j\mid j\leq i\rangle$. Then
 \begin{enumerate}
 \item $L_i\leq L_{i+1}$ for all $i\in \omega$;
 \item $\bigcup_{i\in \omega} L_i=V$;
 \item $L_i+V_i=V$ for all $i\in \omega$.
 \end{enumerate}
 
 Let $B_1$ be a complement to $L_1\cap V_1$ in $L_1$. Then 
 $$B_1+V_1\geq B_1+ (L_1\cap V_1)+V_1=L_1+V_1=V$$
  whence $B_1$ is a supplement to $V_1$. Since $B_1$ is contained in $L_1$ we have $B_1\cap V_1=B_1\cap V_1\cap L_1=0$, that is $B_1$ is a complement to $V_1$ in $V$. 
  Set $\sigma(0)=0$ and $\sigma(1)=1$. The subspace $L_1\cap V_{\sigma(1)}$ has finite dimension so it has trivial intersection with infinitely many of the $V_i$. Define 
  $$\sigma(2)=\min\{i\geq 2\mid L_{\sigma(1)}\cap V_{\sigma(1)}\cap V_{\sigma(2)}=0\}.$$ 

The space $V$ can be written as $V=V_{\sigma(2)}+L_{\sigma(2)}$ and $L_{\sigma(1)}\cap V_{\sigma(2)}=0$. Thus we can choose a complement $B_2$ to $V_{\sigma(2)}\cap L_{\sigma(2)}$ in $L_{\sigma(2)}$, in such a way that $L_{\sigma(1)}\leq B_2 $. In particular $B_1\leq B_2$ and $B_2$ is a complement to $ V_{\sigma(2)}$ in $V$.

Iterating this procedure we end up with a set $\{B_i\mid i\geq 1\}$ and a strictly increasing map $\sigma :\mathbb N\longrightarrow \mathbb N$ satisfying the following conditions
\begin{enumerate}
\item $B_i\leq B_{i+1}$ for al $i\geq 1$;
\item $L_{\sigma(i)}\leq B_{i+1}$ for all $i\geq 1$;
\item $B_i\cap V_{\sigma(i)}=0$ for all $i\geq 1$;
\item $V/V_{\sigma(i)}\simeq B_i$ for all $i\geq 1$.
\end{enumerate}

Finally, put $A_0=0$, $A_1=B_1$ and, when $i>1$, $A_i=B_i\cap V_{\sigma(i-1)}$. For every index $i\geq 1$, we
have $A_i\cap A_{i+1}=A_i\cap B_{i+1}\cap V_{\sigma(i)}\leq B_i\cap  V_{\sigma(i)}=0$  thus $\langle A_i\mid i\geq 1\rangle=\bigoplus_{i\geq 1}A_i$. 
Since $B_2=B_2\cap (B_1+V_{\sigma(1)})$, Dedekind's rule gives $B_2=B_1+(B_2\cap V_{\sigma(1)})=A_1\oplus A_2$. An easy inductive argument shows that $B_n=\bigoplus_{i=1}^nA_i$. 
From this it follows that $L_{\sigma(n)}\leq \langle A_i\mid i\geq 1\rangle$ for all $n$. But $\mathcal B\subseteq \bigcup_{i\in \omega}L_i=\bigcup_{n\in \omega}L_{\sigma(n)}$ whence 
$\langle A_i\mid i\geq 1\rangle=V$. For every given $i\geq 1$ the subspace $V_{\sigma(i-1)}$ can be written as 
$$V_{\sigma(i-1)}=V_{\sigma(i-1)}\cap V=V_{\sigma(i-1)}\cap (V_{\sigma(i)}+B_i)$$
and Dedekind rule can again be invoked to get 
$$V_{\sigma(i-1)}\cap (V_{\sigma(i)}+B_i)=(V_{\sigma(i-1)}\cap B_i)+V_{\sigma(i)}=A_i+V_{\sigma(i)}$$
From this it readily follows that $V_{\sigma(i-1)}/V_{\sigma(i)}\simeq A_i$.
\hfill $\Box$

\

In the next lemma we single out  a rather technical fact to be used in the forthcoming proofs.

\

\begin{lemma}\label{ordered set} Let $\mathcal B$ be a finite set, of order $m$, endowed with a total preorder. Let $\{\mathcal B_i\mid i\in I\}$ be a partition of $\mathcal B$ whose elements have cardinality at most $k$ and  such that the restriction of the preorder  to each $\mathcal B_i$ is  an order. Assume that, for each $i\in I$, we are given an order-preserving injective function $f_i:\mathcal B_i\longrightarrow \Delta=\{1, 2, \dots n\}$ such that
\begin{enumerate}
\item $\Delta=\bigcup_{i\in I}\mathrm{Im}(f_i)$;
\item if $f_i(x)>f_j(y)$ then $x>y$;
\item for each $a\in [2,n]$ there exist $i\in I$ and $x,y\in \mathcal B_i$ such that $f_i(x)=a= f_i(y)+1$.
\end{enumerate}
Then it is possible to find elements $x_l,y_l , \,\, l=1, \dots, r$, with $r= \lfloor (n-2)/k\rfloor$,\- satisfying
\begin{itemize}
\item[(i)] for each $l=1, \dots, r$,  $x_l,y_l\in \mathcal B_{i_l}$ and $\mathcal B_{i_p}\not=\mathcal B_{i_q}$ when $p\not=q$;
\item[(ii)] for each $l=1, \dots, r-1$, $x_{l+1}<y_l<x_l$;
\item[(iii)] for each $l=1, \dots, r$, $f_{i_l}(x_l)=f_{i_l}(y_l)+1$.
\end{itemize}
\end{lemma}

\pf

Choose $i_1\in I$ such that $n,n-1\in \mathrm{Im}(f_{i_1})$ and let $x_1,y_1$ be the preimages of $n$ and $n-1$. Set   $I_2=I\setminus \{i_1\}$ and define  \linebreak $\Delta_2=\Delta\setminus \mathrm{Im}(f_{i_1})$. 
Since $\mathrm{Im}(f_{i_1})$ contains at most $k$ elments, the maximum of $\Delta_2$, $d_2$, is at least $n-k$ and $d_2<n-1$.  Under our assumptions it is possible to find  $i_2\in I_2$, in such a way that $d_2,d_2-1\in \mathrm{Im}(f_{i_2})$. Let $x_2,y_2$  be the preimages, under $f_{i_2}$, of $d_2$ and $d_2-1$.
Since $n-1=f_{i_1}(y_1)>f_{i_2}(x_2)=d_2$, it follows that $y_1>x_2$.
Suppose we have already found a set of  indices $\{i_1, \dots i_s\}\subseteq I$ and elements $x_l,y_l\in \mathcal B_{i_l}$ for $l=1, \dots s$ in such a way that 
\begin{itemize}
\item[(a)] conditions $(i), \, (ii)$ and $(iii)$ hold for each $l=1, \dots s$;
\item [(b)]for all $l=2,\dots,s$ $f_{i_l}(x_l)=d_l=\max \Delta\setminus (\bigcup_{j=1}^{l-1}\mathrm{Im}(f_{i_j}))$ and $f_{i_l}(y_l)=d_l-1$

\end{itemize}

 If $n-sk\geq 2$ let $d_{s+1}$ be the maximum of 
\begin{displaymath}\Delta_{s+1}= (\Delta\setminus ({\bigcup_{l=1}^s}\mathrm{Im}(f_{i_l})))
\end{displaymath}

The set $\Delta_{s+1}$ contains at least $n- ks$ points, so that  $d_{s+1}\geq 2$. Moreover $d_{s+1}<d_{s}-1=f_{i_s}(y_s)$.

By hypothesis it is possible to find  $i_{s+1}$ in such a way that $d_{s+1}, d_{s+1}-1\in \mathrm{Im}(f_{i_{s+1}})$. Let   $x_{s+1},y_{s+1}$ be the  the preimages, under $f_{i_{s+1}}$, of $d_{s+1}$ and $d_{s+1}-1$. It is readily checked that  conditions $(a),\, (b)$ are satysfied for all $i=1, \dots, s, s+1$, so that this procedure can be continued as long as 
  $n-sk\geq 2$,  that is  $s\leq \lfloor (n-2)/k\rfloor$, as claimed.\hfill $\Box$

\begin{lemma}\label{nonunipotent} Let $\s$ be  series of length $n$ in the finite-dimensional $\F$-vector space $V$, and $g$ an element of its stability group $S(\s)$, such that it does not stabilize any proper subseries of $\s$. Assume that $(g-1)^k=0$ for some $k<n-2$.  Then there exists $h\in S(\s)$  such that $(gg^h-1)^{\lfloor (n-2)/k\rfloor-1}\not=0$.
\end{lemma}

\pf The space $V$ can be decomposed as the direct sum $V=\sum_{i=1}^m M_i$ where each $M_i$ is a Jordan block for $g$. Since $(g-1)^k=0$, the dimension of the $M_i$ is bounded by $k$ and this implies that $m\geq \lfloor \dim(V)/k\rfloor\geq \lfloor n/k\rfloor$. For each $i$ we choose a basis $\mathcal B_i=\{v_{i,j}\mid j=1, \dots d(i)\}$ for $M_i$, with respect to which
the matrix representing $g$ is in Jordan canonical form. Thus $v_{i,j}(g-1)=v_{i, j+1}$ if $j<d(i)$, and $0$ otherwise. Let $V_1, V_2, \dots V_n, V_{n+1}$ be the elements of $\s$, with $V=V_1, \, V_{n+1}=0$ and $V_i>V_{i+1}$ for all $i$.   For each vector non-trivial  $v$, there exists a unique $l=l(v)\in \Delta=\{1, 2, \dots , n\}$ such that  $v\in V_l\setminus V_{l+1}$. 

The space  $V$ can be preordered by setting $v\preccurlyeq w$ iff $l(v)\leq l(w)$ and this preorder induces an order on each $\mathcal B_i$. Moreover, if $t$ is in $S(\s)$ and $v\preccurlyeq  w$, then $v(t-1)\preccurlyeq  w(t-1)$.
For $i=1, \dots, r$ let $f_i:\mathcal B_i\longrightarrow \Delta$ be the function defined by $f(v)=l(v)$. Each $f_i$ is order preserving and, if the vectors $x\in \mathcal B_i, \, y\in \mathcal B_j $ satisfy
$ f_i(x)>f_j(y)$, then $x \succ y$. Choose any $a\in [2,n]$. Since $g$ does not stabilize any subseries of $\s$, we have $[V_a,g]\not\subseteq V_{a-2}$, hence there must exist $i$ and $x\in \mathcal B_i$ such that $y=[x,g]\in V_{a-1} \setminus V_{a-2}$. The vector $y$ is in $\mathcal B_i$, thus $f_i(x)=a= f_i(y)+1$. The hypotheses of lemma \ref{ordered set} apply, so this lemma can be invoked to find a sequence $\{x_l, \, y_l \mid l=1, \dots , r=\lfloor (n-2)/k\rfloor\}\subseteq \mathcal B$ satysfying conditions $(i), \, (ii), \, (iii)$ of \ref{ordered set}.
Define a linear transformation $h$ by the following rule: $y_l(h-1)=x_{l+1}$ for all $l=1, \dots , r-1$ and $vh=v$ when $v\in \mathcal B\setminus \{ y_l \mid l=1, \dots , r-1\}$. The transformation $h$ is invertible and belongs to $S(\s)$. It is  important to notice that $(h-1)^2=0$ (hence $h^{-1}=1-(h-1)$ ) and each  $\mathcal B_i$ contains at most one element on which $h$ acts non-trivially. Thus  $h$ acts trivially on each $y_s(g-1)$. The subspace $W=\langle [y_s, g^i]\mid s=1, \dots r, \, i\in \mathbb Z\rangle$ is normalized by $\langle g, h\rangle$, and $\langle g, h\rangle$ acts on $V/W$. The claim will be proved if we show that $(gg^h-1)^r$ is not zero in its action on $V/W$. For each index $s$,  $W\cap M_{i_s}=\langle [y_s, g^i]\mid  i\in \mathbb Z\rangle $ and  $W$ can be decomposed as  $W={\oplus _{s=1}}^r W\cap M_{i_s}$. Whence the set $\{ x_s,y_s\mid s=1, \dots , r\}$ is still linearly independent modulo $W$. Without loss of generality we may  assume, in order to 
 simplify calculations,   that $[y_s,g]=0$ for all $s=1, \dots ,r$.

 Let $s$ be any index between $1$ and $r-2$. We want to understand the action of $gg^h$ on $x_s$ and $y_s$. Using the identity $(gg^h-1)=(g-1)(g-1)^h+(g-1)+(g-1)^h$ we obtain

$$x_s(gg^h-1)=y_s(g-1)^h+y_s+x_s(g-1)h=$$
$$=(y_s-x_{s+1})(g-1)h+y_s+y_s+x_{s+1} 
=2y_s+x_{s+1}-y_{s+1}h=
$$
$$=2y_s+x_{s+1}-y_{s+1}-x_{s+2}$$

and 
$$y_s(gg^h-1)=(y_s-x_{s+1})(g-1)h=-y_{s+1}-x_{s+2}$$
When $s=r-2, r-1$ or $r$ the above calculations can be easily modified. For $s=1, \dots , r$ define $T_s=\langle x_i, y_i\mid i\geq s\rangle$ and $B_s=\langle y_s, x_i, y_i\mid i> s\rangle$. The above calculations show that $gg^h$ stabilizes the series 
$$\mathcal S: 0<B_r<T_r<\dots  <B_1<T_1$$
In particular $[T_s,gg^h]\leq B_s$ and $[B_s, gg^h]\leq B_{s+1}$. An easy inductive argument shows that $y_1(gg^h-1)^l=(-1)^ly_{l+1}+v_l$ where $v_l$ is a suitable element of $T_{l+2}$. In particular  $y_1(gg^h-1)^{r-1}\not=0$, and the claim is proved. \hfill $\Box$

\bigskip 

Lemma \ref{nonunipotent} must be extended in order to cover the case of vector spaces of infinite dimension.

\

\begin{lemma}\label{nonunipotent2} Let $\s$ be  a series  in the  $\F$-vector space $V$, $k,n\geq 1$ and $g$ an element of the  stability group $S(\s)$, such that  $(g-1)^k=0$. If  $g$ does not stabilize any  subseries of $\s$ of length smaller than $n$ and $k<n-2$, there exists $h\in S(\s)$  such that $(gg^h-1)^{\lfloor (n-2)/k\rfloor-1}\not=0$.
\end{lemma}

\pf Let $\s_0:0<V_n<V_{n-1}<\dots <V_1=V$ be any subseries of $\s$ of length $n$. Our hypothesis imply that, for each $i=1, \dots, n-1$, there exists $v_i\in V_i$ such that $[v_i,g]\in V_{i+1}\setminus V_{i+2}$. Consider the subspace $W_0=\langle v_1, v_2, \dots v_n\rangle$ and construct $W=W_0+\sum_{i}W_0(g-1)^i$. The subspace $W$ is $\langle g\rangle$-invariant and has finite dimension. The series $\s_0$ induces a series $\mathcal S$ in $W$. This series has length $n$, and  $g$ does not stabilize any proper   subseries of $\mathcal S$. Lemma \ref{nonunipotent} can be invoked to find an element $t$ in the stability group of $\mathcal S$, such that $(gg^t-1)^{\lfloor (n-2)/k\rfloor-1}\not=0$ when acting on $W$. Since $W$ is finite dimensional, it is an easy matter to extend $t$ to an element $h\in S(\s)$. Clearly $(gg^h-1)^{\lfloor (n-2)/k\rfloor-1}\not=0$ and the lemma  is proved. \hfill $\Box$

\bigskip

We are now able to describe the Hirsch-Plotkin radical for stability groups of ascending series.

\

\begin{thm}\label{ascending}Let $\s$ be a series in $V$, $S(\s)$ its stability group and $H(\s)$ the Hirsch-Plotkin radical of $S(\s)$. If the order type of $\s$ is an ordinal, then $H(\s)=\{ g\in S(\mathcal L) \mid g\textrm{ stabilizes a finite subseries of } \mathcal L\}$.
\end{thm}

\pf By way of contradiction assume there exists $g$ in $H(\s)$, such that it does not stabilize any finite subseries of $\s$. 
The series $\s$ can be written as $\s=\{V_\alpha\mid \alpha<\lambda\}$ for a suitable ordinal $\lambda$. The set $C_1=\{\alpha \mid [V_\alpha,g]=0\}$ contains at least the ordinal $1$, so $W_1=\sup\{V_\alpha\mid \alpha\in C_1\}$ is a non-trivial element of $\s$. This process can be applied again to $V/W_1$  in order to produce, using an inductive argument, a subset $\s_0=\{W_i\mid i\in \omega\}\subseteq \s$. This subset must be infinite, because $g$ does not stabilize any finite subseries of $\s$. Set  $W=\cup_{i\in \omega}W_i$ and consider $g$ in its action on $W$. For any given $m<n$, the series $\s_0$ induces a series $\s(n,m)$ of length $n-m$ in $W_n/W_m$. The element $g$ stabilizes $\s(n,m)$ in its action on $W_n/W_m$, but it is clear from the definition of $\s_0$, that no proper subseries of $\s(n,m)$ can be stabilized by $g$. We now define the series $\s_1\subseteq \s_0$ setting $U_1=W_1$ and, when $i>1$, defining $U_i$ to be the unique element of $\s_0$ such that the 
 series induced by $\s_0$ on $U_i/U_{i-1}$ has length $i$. Choose $A_i$ a complement to $U_{i-1}$ in $U_i$ and fix any isomorphism $\sigma_i: U_i/U_{i-1}\longrightarrow A_i$. These isomorphisms can be used to define an action of $g$ on each $A_i$. Morever the subspaces $A_i$  can be endowed with a series $\mathcal S_i$, which is the image, via $\sigma_i$, of the series induced by $\s_0$ on $U_i/U_{i-1}$. In particular $\left|\mathcal S_i\right|=i$. Call $g_i$ the transformation induced by $g$ on $A_i$. Each $g_i$ is in the stability group of $\mathcal S_i$ and $(g_i-1)^k=0$. Use lemma \ref{nonunipotent2} to find, for each $i>k$, an element $h_i$ in the stability group of $\mathcal S_i$, such that $(g_i{g_i}^{h_i}-1)^{\lfloor i/k\rfloor-1}\not=0$ on $A_i$. Since  $W=U_1\oplus(\bigoplus_{i>1}A_i)$ the space $V$ can be written as 
$V=M\oplus U_1\oplus(\bigoplus_{i>1}A_i)$ once a suitable complement $M$ has been selected. It is then possible to define an automorphism $h$ of $V$ setting $h=1$ on $M\oplus U_1\oplus(\bigoplus_{i\leq k}A_i)$ and $h=h_i$ on $A_i$ when $i>k$. The element thus defined  lies in the stability group of $\s$ and $gg^h$ is not unipotent because, for each $m\in \mathbb N$, there is a suitable section $U_i/U_{i-1}$ on which $ (gg^h-1)^m$ acts non-trivially. But $g^h\in H(\s)$, hence $gg^h$ should be unipotent. This contradiction proves that $g$ must stabilize a finite subseries of $\s$, showing that the  claim holds. \hfill $\Box$

\bigskip

In the next theorem we drop the assumption on the order type of the series  but, on the other hand, we need to restrict ourselves to vector spaces of countable dimension.

\

\begin{thm}\label{countable}Let $\s$ be a series in $V$, $S(\s)$ its stability group and $H(\s)$ the Hirsch-Plotkin radical of $S(\s)$. If the dimension of $V$ is at most $\aleph_0$, then $H(\s)=\{ g\in S(\mathcal L) \mid g\textrm{ stabilizes a finite subseries of } \mathcal L\}$.
\end{thm}

\pf We argue by contradiction, assuming there exists $g$ in $H(\s)$, such that it does not stabilize any finite subseries of $\s$. 
If $U,W$ are elements of $\s$ with $U<W$, the group $S(\s)$ acts on $W/U$ as the full stabilizer of $\s\cap W/U$ so that, in order to get the desidered contradiction, it will be enough, if necessary, to consider the action of $g$ on a suitable section of $\s$.

    Set $W_0=V$ and $W_1=\min \{V_\alpha\mid V_\alpha\in \s, \, g \textrm{ acts trivially on } V/V_\alpha \}$. Once $W_n$ has been defined, set 
$$W_{n+1}=\min \{V_\alpha\mid V_\alpha\in \s, \, g \textrm{ acts trivially on } W_n/V_\alpha \}.$$ Two cases should be considered separately.

\

{\bf Case 1.} \emph{ The set $\mathcal W=\{W_i\mid i\in \omega\}$ is finite.}

\

By restricting to the lowest term of $\{W_i\mid i\in \omega\}$, we may assume that this set contains only $V$.
Let $U_1=0$ and $U_2$ be  any proper element of $\s$, such that $[U_2,g]\not=0$. This subspace does exist because $g$ is not the identity. Moreover $g$ does not stabilize any subseries of $\s\cap U_2$ of length smaller than $2$.
Assume we have already found $U_1, U_2, \dots U_n\in \s$ such that,
\begin{itemize}
\item $U_{i-1}<U_i$ for all $i=2, \dots, n$ and $U_n<V$;
\item for each $i=2, \dots ,n$, $g$ does not stabilize any subseries of $\s\cap U_i/U_{i-1}$ of length smaller than $i$.
\end{itemize} 
 The transformation $g$ does not stabilize any finite subseries of $\s/U_n$ 
in its action on $V/U_n$, otherwise $\mathcal W$ would contain more than one element. Hence we define $U_{n+1}$ to be  any proper element of $\s$, such that $g$ does not stabilize any finite subseries of $\s\cap U_{n+1}/U_n$ of length smaller than $n+1$. This process gives an ascending subseries $\mathcal U=\{U_i\mid i\geq 1\}\subseteq \s$ such that, for every $i\geq 2$,  $g$ does not stabilize any subseries of $\s\cap U_i/U_{i-1}$ of length smaller than $i$.
It is then possible to apply the same argument used in \ref{ascending}, to come up with an $h\in S(\s)$ such that $gg^h$ is not unipotent. This proves that case 1 can not occur.

\

{\bf Case 2.} \emph{ The set $\mathcal W=\{W_i\mid i\in \omega\}$ is infinite.}

\

The subspace $W=\cap_{i\in \omega}W_i$ belongs to $\s$, and  $S(\s)$ acts on $V/W$ as the stability group of the series $\s/W$. We don't loose generality assuming $W=0$. Lemma \ref{decomposition} can be used to find a strictly increasing function $\sigma:\omega\longrightarrow \omega$ and subspaces $A_i$, $ i\in \omega$,  such that 
\begin{enumerate}
\item $V=\bigoplus_{i\geq 1}A_i$;
\item $A_i\simeq W_{\sigma(i-1)}/W_{\sigma(i)}$ for all $i\geq 1$.
\end{enumerate}
 Let $n_i$ be the length of the series induced by $\s_0=\{W_i\mid i\in \omega\}\cup \{0\}$ on $W_{\sigma(i-1)}/W_{\sigma(i)}$. It is clear that we can choose $\sigma$ in such a way that the sequence $\{n_i \mid i\in \omega\}$ is unbounded. Using the same technique used in theorem \ref{ascending},  we find $h\in S(\s)$ such that $gg^h$ is not unipotent, a contradiction.

This contradiction shows that the claim holds.\hfill $\Box$

\bigskip

From the proof of theorem \ref{countable} it should be clear that, in order to drop the assumption on the dimension of $V$, we should be able to prove our claim when $\s$ is a descending chain. A typical example occurs when $V$ is a subspace of $\F^\omega$ containing $\mathrm{Dir}_{i\in \omega}\F$, and the series $\s=\{V_i\mid i\in \omega\}\cup \{0\}$ is defined by  $V_i=V\cap \{ (a_j)\mid a_j=0 , \forall j<i\}$. When $V$ has uncountable dimension, we have been unable to adapt our technique to this setting. However our machinery can be used when a particular kind of basis exists.

Let $\mathcal B$ be a basis of $V$ and, for each jump $j=(B,T)$ of $\s$ define $\mathcal B_j=\{v\in \mathcal B \mid v \textrm{ belongs to } j\}$. The basis $\mathcal B$ is said to be \emph{$\s$-adapted} if, for  each jump $j=(B,T)$, the set $\mathcal B_j+B=\{v+B\mid v\in \mathcal B_j\}$ is a basis of $T/B$. 
It is worth remarking that the existence of $\s$-adapted basis, is not granted in infinite dimension. The easiest example is perhaps the following.

\

Choose any field $\F$ and let $V=\F^\omega$. For each $i\in \omega $ set $V_i=\{(v_j)\mid v_j=0 \, \, \forall j<i\}$. The series $\s=\{V_i\mid i\in \omega\}\cup \{0\}$ has countably many jumps, each of dimension one, while $V$ has uncountable dimension. Hence no $\s$-adapted basis exists.

\

We point out an important property of $\s$-adapted basis

\bigskip

\begin{lemma}\label{adapted section} Assume that $\mathcal B$ is an $\s$-adapted basis of $V$, and $U<W$ two elements of $\s$. Then the set $\mathcal B(W/U)=\{v+U\mid v\in (W\setminus U)\cap \mathcal B\}$ is a basis for $W/U$.
\end{lemma}

\pf Let $M=\langle (W\setminus U)\cap \mathcal B \rangle $. If $M+U<W$ choose  $w\in W\setminus (M+U)$ in such a way that, if $w=\sum_{v\in \mathcal B}\lambda_vv$, the support of $w$ $S(w)=\{v\mid \lambda_v\not=0\}$ has minimal cardinality. Pick $s\in S(w)$ and consider $w_0=w-\lambda_ss$. The jump $(B,T)$ to which $s$ belongs, clearly satysfies $U\leq B<T\leq W$, hence $w_0$, if not $0$,  is still contained  in $W\setminus M+U$. Since $S(w_0)$ is strictly contained in $S(w)$, $w_0=0$ showing that $w\in M$, a contradiction.\hfill $\Box$

\bigskip

\begin{lemma}\label{extending} Assume that $\mathcal B$ is an $\s$-adapted basis of $V$ and $\mathcal F=\{ W_\lambda/U_\lambda\mid \lambda\in \Lambda\}$ a set of sections of $\s$ such that, for any pair of distinct indices $\alpha,\beta\in \Lambda$, we have $W_\beta\leq U_\alpha$ or 
$W_\alpha\leq U_\beta$. For each $\lambda\in \Lambda$ let $\s_\lambda$ be the series 
$W_\lambda\cap \s/U_\lambda$. Assume we are given, for each $\lambda\in \Lambda$, an element  
$h_\lambda\in S(\s_\lambda)$. Then there exists $h\in S(\s)$ such that, $\forall \, \lambda\in \Lambda$,   $h$ induces $h_\lambda$ in its action on $W_\lambda/U_\lambda$.
\end{lemma}

\pf To describe  $h$ it is enough to define it on the members of $\mathcal B$. 
Let $\mathcal B(\lambda)=\{v\in \mathcal B\mid v\in W_\lambda\setminus U_\lambda\}$. If $v\not\in \bigcup_{\lambda\in \Lambda}\mathcal B(\lambda)$, set $vh=v$. If $v\in \mathcal B(\lambda)$ and $(v+U_\lambda)h_\lambda=\sum_{u\in \mathcal B(\lambda)}a_u(u+U_\lambda)$, set $vh=
\sum_{u\in \mathcal B(\lambda)}a_uu$.
The map $h$ belongs to $\gl{V}{\F}$. If $vh=0$ write $v=\sum_{\lambda\in I}v_\lambda+w$ where $I$ is a suitable finite subset of $\Lambda$,  each $v_\lambda$ has support contained in  $\mathcal B(\lambda)$ and $w$ has support in 
$\mathcal B\setminus \bigcup_{\lambda\in \Lambda}\mathcal B(\lambda)$. The equation 
$$0=vh=\sum_{\lambda\in I}v_\lambda h+wh=\sum_{\lambda\in I}v_\lambda h+w$$
 holds if and only if $w=0$ and $v_\lambda h=0$ for all $\lambda \in I$. Thus $v_\lambda=0\, \, \forall \, \lambda\in I$, showing that $h$ is injective.
To prove surjectivity select $v\in V$ and, using the same notation of the above paragraph, write $v=\sum_{\lambda\in I}v_\lambda+w$. For each $\lambda \in I$ there exists $u_\lambda\in \langle \mathcal B(\lambda)\rangle$, such that $(u_\lambda+U_\lambda)h_\lambda=v_\lambda+U_\lambda$, hence $u_\lambda h=v_\lambda$. Since $wh=w$ we have $(\sum_{\lambda\in I}u_\lambda+w)h=v$.
The fact that $h$ is in $S(\s)$ is clear.
\hfill $\Box$

\bigskip

\bigskip
In presence of  $\s$-adapted basis, we are able to prove our conjecture.

\bigskip

\begin{prop}\label{adapted} Let $\s$ be a series in $V$, $S(\s)$ its stability group and $H(\s)$ the Hirsch-Plotkin radical of $S(\s)$. If there exists an $\s$-adapted basis, then $H(\s)=\{ g\in S(\mathcal L) \mid g\textrm{ stabilizes a finite subseries of } \mathcal L\}$.
\end{prop}

\pf  By way of contradiction assume that $g$ belongs to $H(\s)$ but it does not stabilize any finite subseries of $\s$, and let $k$ be such that $(g-1)^k=0$. We use the same approach of   theorem \ref{countable}, and the notation established therein. Case 1 can be handled without any change.

If we are in case 2, we choose a subset $\{U_i\mid i\in \omega\}\subseteq \mathcal W$, in such a
 way that $g$, in its action on $U_i/U_{i+1}$,  does not stabilize any subseries of $\s_i=\s\cap 
 U_i/U_{i+1}$, of length less than $i$. If $g_i$ stands for the transformation induced by $g$
  on $U_i/U_{i+1}$, lemma \ref{nonunipotent2} can be used to produce elements $h_i\in S(\s_i)$,
   such that $(g_i{g_i}^{h_i}-1)^{\lfloor i/k\rfloor-1}\not=0$. We invoke now lemma \ref{extending} to find $h\in S(\s)$ whose restriction to each $U_i/U_{i+1}$ is $h_i$. Thus $gg^h$ is not unipotent and this contradiction gives the claim.\hfill $\Box$

\bigskip   
   
What we have proved about $\s$-adapted basis, can be used to discuss a quite different case. The  situation we want to investigate is best described in terms of topology.

\

Let $V$ be an $\F$-vector space, and $\s$ a series in $V$. There exists a unique topology $\tau_\s$ on $V$, admitting  $\s\setminus \{0\}$ as a basis for the open neighborhoods of $0$. We call this the \emph{ $\s$-topology} on $V$. Clearly $\tau_\s$ is $T_2$ if and only if $0=\bigcap\{ U\mid U\in \s \, \,\textrm{ and } U\not= 0\}$. 

Recall that, in this situation, we have a notion of \emph{Cauchy net}, so that it is meaningful to talk about completeness of the topological space $(V, \tau_\s)$.

\

\begin{lemma}\label{extension}Let $V$ be an $\F$-vector space, $\s$ a series in $V$ and assume that $(V,\tau_\s)$ is complete and $T_2$. If $W$ is a dense subspace and $h\in \gl{W}{\F}$ belongs to $S(W\cap \s)$, there exists a unique $\overline{h}\in S(\s)$ such that $wh=w\overline{h}$ for all $w\in W$.
\end{lemma}

\pf The map $h$ is continous with respect to the topology induced by $\s\cap W$ on $W$, because 
$h$ normalizes every element of $\s\cap W$. If $X$ is 
any subset of $V$, its closure is $\overline{X}=\cap \{X+U\mid U\in \s\, \textrm{ and } U
\not=0\}$. Hence $W+U=V$ for all $0\not=U\in\s$. Let $v\in V$ be any vector and choose any net $
\{ v_\lambda\mid \lambda\in \Lambda\}$ in $W$,  converging to $v$. The net $\{v_\lambda h\mid \lambda\in 
\Lambda\}$ is Cauchy. In fact choose $0\not=U\in \s$ and let $\lambda\in \Lambda$ be such that,
 if $\alpha, \beta\geq \lambda$, then $v_\alpha-v_\beta\in U\cap W$.  Thus $v_\alpha h-v_\beta 
 h=(v_\alpha-v_\beta)h$ belongs to $(U\cap W)h=U\cap W$. This shows that the net is Cauchy in $
 (V,\tau_\s)$, so that it has a unique limit $v\overline {h}=\lim_{\lambda\in \Lambda} v_\lambda h$. 
 It is readily seen that $v\overline {h}$ does not depend on the choice of the approximating 
 net $\{v_\lambda\mid \lambda\in \Lambda\}$, so that the function 
 $\overline{ h}:V\longrightarrow V$ is well defined. Suppose  $\ker (\overline{h})\not=0$ and choose a non trivial 
 $v\in \ker (\overline{h})$. There exists a jump $(B,T)$ of $\s$ such that $v\in T\setminus B$.
  The density of $W$ allows to find  $w\in W$ such that $w+B=v+B$. 
  Hence $B=v\overline{h}+B=wh+B$, and $wh\in B\cap W$. This is impossible because $h$ acts trivially on $T\cap W/B\cap W$, whence  $\overline{h}$ is injective.
  To prove surjectivity write $v$ as the limit of the Cauchy net $\{v_\lambda \mid \lambda\in 
\Lambda\}$ in $W$, and notice that $\{u_\lambda=v_\lambda h^{-1}\mid \lambda\in 
\Lambda\}$ is still a Cauchy net. If $u$ is its limit we get $u\overline{h}=v$.
Consider the jump $(B,T)$ of $\s$. The density of $W$ gives $T=B+(T\cap W)$, so that 
$$[T,\overline{h}]=[B,\overline{h}]+[T\cap W,\overline{h}]=[B,\overline{h}]+[T\cap W,h]\leq B
$$
Hence $\overline{h}$ belongs to $S(\s)$. The uniqueness of $\overline{h}$ is clear.
\hfill $\Box$

\

We are now ready to prove our conjecture for complete spaces.
\bigskip

\begin{thm}\label{complete} Let $V$ be an $\F$-vector space, $\s$ a series in $V$ and assume that $(V, \tau_\s)$ is a complete topological space. Then $H(\s)=\{ g\in S(\mathcal L) \mid g\textrm{ stabilizes a finite subseries of } \mathcal L\}$.
\end{thm}

\pf For each jump $j=(B,T)$ select a subset $\mathcal B_j$ of $V$ such that $\{v+B\mid v\in \mathcal B_j\}$ is a base for $T/B$. The subspace $W=\langle \mathcal B_j \mid j \textrm{ is a jump for } \s\rangle$ is dense in $(V, \tau_\s)$. In fact, if  this is  false, it is possible to find   a non-trivial $K\in \s$ for which  $W+K<V$. Choose  $v\in V\setminus (W+K)$. The vector $v$ belongs to $T\setminus B$ for a suitable jump $j=(B,T)$, and $K$ is clearly strictly contained in $T$. Since $\s$ is linearly ordered, either $B<K$ or $K\leq B$. But $(B,T)$ is a jump, so there are no elements of $\s$ between $B$ and $T$, thus $K\leq B$. We can write $v=w+k$  for suitable vectors $w\in W$ and $k\in K$, because $W$ contains $\mathcal B_j$, so that $v\in W+K\leq W+U$, a contradiction. 
If $R/S$ is any section of $\s$, the usual application of Dedekind's rule gives  $R=R\cap (W+S)=(R\cap W)+S$. 
Hence $R/S\simeq R\cap W/S \cap W$.
Let $g$ be any element of $H(\s)$ and, by way of contradiction, assume it does not stabilize any finite subseries of $\s$. We follow the argument of theorem \ref{countable} and stick to the notation there defined. The discussion of case 1 remains unchanged. 

If we are in case 2 we choose a subset $\{U_i\mid i\in \omega\}\subseteq \mathcal W$, in such a
 way that $g$, in its action on $U_i/U_{i+1}$,  does not stabilize any subseries of $\s_i=\s\cap 
 U_i/U_{i+1}$, of length less than $i$. If $g_i$ stands for the transformation induced by $g$
  on $U_i/U_{i+1}$, lemma \ref{nonunipotent2} can be used to produce, when $i-2>k$, elements $h_i\in S(\s_i)$,
   such that $(g_i{g_i}^{h_i}-1)^{\lfloor (i-2)/k\rfloor-1}\not=0$. For each $i>k+2$, 
   $U_i/U_{i+1}$ is naturally isomorphic to $ U_i\cap W/U_{i+1}\cap W$, and each $h_i$ induces, 
   via this isomorphism, a   transformation $k_i$ in the stabilizer of the series induced by $\s$ 
   on the  section  $ U_i\cap W/U_{i+1}\cap W$.
   We invoke now lemma \ref{extending} to find $k\in S(\s\cap W)$ whose restriction to each 
   $ U_i\cap W/U_{i+1}\cap W$ is $k_i$. By lemma \ref{extension} $k$ has a unique extension $h$ to the whole $V$, and $h\in S(\s)$. As  it is readily seen, $h$ induces $h_i$ in its action on $U_i/U_{i+1}$ so that, by the usual argument, $gg^h$ is not unipotent, contradicting the fact that $g$ belongs to $H(\s)$. This contradiction proves that the claim holds.
   \hfill $\Box$

   \bigskip

\section{ Fitting groups}
\bigskip

In this section we discuss a particular kind of unipotent representations for Fitting groups. The results obtained are then used to prove the following theorem

\

\begin{thm}\label{fitting as hp} Let $G$ be a Fitting group such that $G=X^G$ for some countable subgroup $X\leq G$. If $G$ is torsion-free or a $p$-group for some prime $p$, then there exist a field $\mathbb F$, an $\mathbb F$-vector space $V$ and a series $\mathcal L$ in $V$, such that $G$ can be embedded in the Hirsch-Plotkin radical of $S(\mathcal L)$.
\end{thm}

\

In particular the above theorem holds when $G$ is countable.
 
\begin{Def} A \emph{preorder} on a set $\Omega$ is a reflexive and transitive binary relation. \end{Def}

We need a week notion of maximality  in preordered set.

\begin{Def} Let $(\Omega, \leq)$ be a preordered set. An element $m\in \Omega$ is said to be \emph{nearly maximal} if, whenever $m\leq a$, we also have $a\leq m$.\end{Def}

When $(\Omega, \leq )$ is a preordered set, define  $a\sim b$ if and only if $a\leq b$ and $b\leq a$. This is clearly an equivalence relation. Let  $\Delta=\Omega/\sim$ be the quotient set and choose $[a], [b]\in \Delta$. If $a\leq b$ then $x\leq y$ whenever $x\in [a]$ and $y\in [b]$. Therefore  setting $[a]\preccurlyeq [b] \, \Longleftrightarrow \, a\leq b$ defines an order on $\Delta$. We shall refer to $(\Delta, \preccurlyeq)$ as to the \emph{ canonical ordered set} associated to $(\Omega, \leq)$.

\begin{thm}\label{rep} Let $G$ be a Fitting group. If $G$ is torsion-free or a $p$-group for some prime $p$, then there exist a field $\mathbb F$, an $\mathbb F$-vector space $V$ and a series $\mathcal L$ in $V$, such that $G$ can be embedded in $S(\mathcal L)$ in such a way that  for each $g\in G$  there exists $n=n(g)\in\mathbb N$, such that $[V, _n g^G]=0$.
\end{thm}

Clearly the requirement that $G$ is torsion-free or a $p$-group can not be dropped, since the Hirsch-Plotkin radical of $S(\mathcal L)$ is torsion-free or a $p$-group according to $\mathbb F$ being of characteristic $0$ or $p$.
We need some  lemmata.

\begin{lemma}\label{zorn} Let $(\Omega,\leq)$  be a preordered set, such that every chain in $\Omega$ has an upper bound. Then $(\Omega,\leq)$ has nearly maximal elements.\end{lemma}

\pf Let $(\Delta,\preccurlyeq)$ be the canonical ordered set associated to $(\Omega, \leq)$.
Choose any chain $\mathcal L=\{\lambda_i\mid i\in I\}\subseteq \Delta$ and, for each $i\in I$, choose $l_i\in \lambda_i$. Then $L=\{ l_i\mid i\in I\}$ is a chain in $\Omega$, so that there exists an upper bound $l$. Thus the element $\lambda=[l]$ is an upper bound for $\mathcal L$ in $\Delta$. It is then possible to apply Zorn's lemma to the ordered set  $(\Delta, \preccurlyeq)$ proving that it possesses at least one maximal element $\mu$. It is readily seen that any $m\in \mu$ is a nearly maximal element in $(\Omega, \leq)$.\hfill $\Box$

\bigskip

\begin{lemma}\label{fitting} Let $G\leq\gl{V}{\mathbb F}$ be such that, for each $g\in G$, there exists $n=n(g)\in \mathbb N$ satysfying $[V, _n g^G]=0$. Then there exists a series $\mathcal L$ in $V$ with $G$ contained in  $S(\mathcal L)$.
\end{lemma}

\pf The group $G$ normalizes the series $\{0, V\}$, hence we choose a series $\mathcal L$ maximal subject to being normalized by $G$. If $G$ is not contained in $S(\mathcal L)$, there exists a jump $(B,T)$ such that $[T,G]$ is not contained in $B$. In particular there must exist $g\in G$ such that $[T,g]\not\leq B$ and, therefore, $[ T, g^G]\not\leq B$. On the other hand, if $n=n(g)$, we have $[T,_n g^G]\leq [V, _n g^G]=0$ so that $g^G$ acts non-trivally on $T/B$ and $B<B+[T,g^G]<T$, because $n>1$. Since $G$ normalizes $B+[T,g^G]$, the group $G$ normalizes the series $\mathcal L\cup \{ B+[T,g^G]\}$, contradicting the fact that $\mathcal L$ is maximal among the series normalized by $G$. 
\hfill$\Box$

\bigskip

We can now prove theorem \ref{rep}. 

\pf Since the theorem is true when $G$ is finite, we assume henceforth that $\left|G\right| \geq \aleph_0$. 

Let $N$ be a normal subgroup of $G$, $\mathbb F$ a field and $V$ an $\mathbb F$- vector space. We say that an injective homomorphism $\sigma: N\longrightarrow \gl{V}{\mathbb F}$ is an $\mathcal F$-representation for $N$ if, for every $x\in N$, there exists $n=n(x)\in \mathbb N$, such that $[V, _n (x^G)\sigma]=0$. Notice that, since $x^G$ acts unipotently on $V$, the field $\mathbb F$ has characteristic $0$ or $p$ according to $G$ being torsion-free or a $p$-group. Let $\mathbb F$ be $\mathbb Q$ or the field with $p$ element, according to $G$ being torsion-free or a $p$-group, and let $V$ be an $\mathbb F$-vector space of dimension $\left|G\right|$. Consider the set 
$$\mathcal N=\{ (N, \sigma_N)\mid N\trianglelefteq G \textrm{ and } \sigma_N:N\rightarrow \gl{V}{\mathbb F}
\textrm{ is an } \mathcal F-\textrm{representation }\}
$$

Of course $\mathcal N$ is not empty and we shall define a preorder on it. For $(N,\sigma_N), \, (M, \sigma_M)\in \mathcal N$ we say that 
$(N,\sigma_N),\preccurlyeq (M, \sigma_M)$ if   $N\leq M$ and, whenever $g\in N$ and $[V, _n (g^G)\sigma_N]=0$, then $[V, _n (g^G)\sigma_M]=0$ too.
Our next aim is to show that $(\mathcal N, \preccurlyeq)$ has nearly maximal elements. In view of lemma \ref{zorn} it is enough to show that every chain in $(\mathcal N, \preccurlyeq)$ has an upper bound. 
Let $\mathcal C=\{(N_i, \sigma_i)\mid i\in I\}$ be a chain in $\mathcal N$. For each $i\in I$ define $A(i)=\{j\in I\mid N_i\leq N_j\}$. Each element of $\mathcal A=\{ A(i)\mid i\in I\}$ is  non-empty and the intersection of finitely many of them still belongs to $\mathcal A$. Thus $\mathcal A$ is contained in an ultrafilter $\mathcal U$ on $I$.  The ultrapower  $W=V^I/\mathcal U$ is a vector space over the field $\mathbb K= {\mathbb F}^I/\mathcal U$ and the ultraproduct $\Cr\{\gl{V}{ \mathbb F}\mid i\in I\}/\mathcal U$ can be seen as a subgroup of $\gl{W}{\mathbb K}$ in a natural way.  
If $\mathbb F$ is the field with $p$ elements, then $\mathbb F\simeq \mathbb K$. In fact the first-order sentence $``\forall a \rightarrow \, a^p=a"$ holds in every component of ${\mathbb F}^I$ and so  it holds in the ultrapower $\mathbb K$. Hence every element of $\mathbb K$ is a root of the polynomial $x^p-x$, thus showing that $\mathbb K$ is isomorphic to $\mathbb F$.  When $\mathbb F=\mathbb Q$, the field  $\mathbb K$ has characteristic $0$, hence it  contains a subfield isomorphic to $\mathbb Q$. In both cases we have that $\mathbb K$ contains a subfield isomorphic to $\mathbb F$, and  we shall view $W$ as an $\mathbb F$-vector space.
Moreover the group $\Cr\{\gl{V}{ \mathbb F}\mid i\in I\}/\mathcal U$ will be seen as a subgroup of  $\gl{W}{\mathbb F}$.

Let $N=\cup_{i\in I}N_i$. This is a normal subgroup of $G$. Given $x\in N$ and $i\in I$ set $(x)\tau_i=1$ if $x\not\in N_i$ and $(x)\tau_i=(x)\sigma_i$ otherwise. Consider the map $\tau :N\longrightarrow \Cr\{\gl{V}{ \mathbb F}\mid i\in I\}$ defined as $(x)\tau=((x)\tau_i)_{i\in I}$ and set $\eta=\tau \pi$, where 
$$\pi:\Cr\{\gl{V}{\mathbb F}\mid i\in I\}\longrightarrow \Cr\{\gl{V}{ \mathbb F}\mid i\in I\}/\mathcal U\leq \gl{W}{\mathbb F}$$
 is the canonical projection. It is an easy matter to show that $\eta$ is an injective homomorphism, so that it embeds $N$ into $\gl{W}{\mathbb F}$. Choose any $x\in N$. If $x\in N_i$ there exists $n\in \mathbb N$ such that $[V, _n (x^G)\sigma_i]=0$. Since $\mathcal C$ is a chain in $\mathcal N$, we have $[V, _n (x^G)\sigma_j]=0$ for all $j\in A(i)$, hence $[W,_n (x^G)\eta]=0$. Thus $\eta$ is an $\mathcal F$-representation for $N$.
For each $1\not=g\in G$ choose $v_g\in W$ in such a way that $[v_g,(g)\eta]\not=0$, and let $M$ be the $\mathbb FG$-submodule of $W$ generated by $\{ v_g\mid g\in G\}$. Since $\dim(M)\leq \left|G\right|$, $M$ can be embedded into  $V$. Hence $\eta$ affords an $\mathcal F$-representation $\sigma: N\longrightarrow V$ and it is readily seen that $(N,\sigma)$ is an upper bound for $\mathcal C$.
 Lemma \ref{zorn} can now be invoked to produce a nearly maximal element  $(H, \tau)$ in $\mathcal N$.

  To prove the theorem it is enough to show that $H=G$ and we shall prove this equality arguing by contradiction. If $\tau$ embeds $H$ into $\gl{V}{\F}$, we shall consider, from now onward, $H$ as a subgroup of $\gl{V}{\F}$, in order to simplify notation.
If $H$ is a proper subgroup $G$, there exists a nilpotent normal subgroup $N\trianglelefteq G$ such that $H<HN$. To see this   choose $x\in G\setminus H$ and set $N=x^G$. The subgroup $N$ is nilpotent because $G$ is a Fitting group. Set $K=HN$ and $L=H\cap N$. Let $\omega$ be the augmentation ideal of $\F N$ and, for $k\in \mathbb N$ consider $\F N/\omega^k$. We choose $k$ in such a way that $\F N/\omega^k$ is a faithful $N$-bimodule. This is indeed possible, see \cite{hartley} theorems 2.2.1 ans 2.2.2. Since $\F N/\omega^k$  is an $L$-bimodule,  we can  construct the group $M=V\otimes_L (\F N/\omega^k)$ which can be endowed with a structure of right $L$-module in the usual way, setting $v\otimes \alpha.x=v\otimes (\alpha x)$ for all $x\in L$, $v\in V$, and $\alpha\in \F N/\omega^k$, and extending by linearity. Notice that the dimension of $M$ as a vector space over $\mathbb F$ is at most $\left| G\right |$.

We prove, as a first step, that $M$ is a faithful $L$-module. To this extent we notice that $\F N$ is a free $L$-module and, if $T\cup\{1\}$ is a transversal for $L$ in $N$ , the set $\{ 1, x-1\mid x\in T\}$ is an $L$-basis for $\F N$. Thus every element in $V\otimes_L \F N$ can be uniquely written as   
$(v_1\otimes 1)+\sum_{x\in T}v_x\otimes (x-1)$ (see \cite{hung} chap. IV theorem 5.11). 

Starting with the exact sequence of  $L$-bimodules
$$
0\longrightarrow \omega^k\longrightarrow \F N \longrightarrow \F N/\omega^k \longrightarrow 0
$$
and tensoring over $L$ with the left $L$-module $V$, we get the exact sequence of abelian groups

$$
V\otimes_L\omega^k\stackrel{\sigma}
{\longrightarrow} V\otimes_L \F N \longrightarrow M\longrightarrow 0
$$
It is readily seen that the homomorphisms in the last  sequence are homomorphisms in the category Mod-$L$. We stick our attention to the map $\sigma$. We have $(\sum_{i}v_i\otimes \alpha_i)\sigma=\sum_{i}v_i\otimes \alpha_i$. Since $\alpha_i\in \omega^k$ for all $i$, we have that 
$$\mathrm{Im}(\sigma)\leq V\otimes_L \omega=\oplus_{x\in T}V\otimes_L(x-1).$$ In order to show that $M$ is faithful, we use the isomorphism \linebreak $M\simeq (V\otimes_L\F N)/\mathrm{Im}(\sigma)$. Let $g\in L$ be an element acting trivially on $ V\otimes_L\F N/\mathrm{Im}(\sigma)$. Hence, for every $v\otimes \alpha\in V\otimes_L \F N$, we have $v\otimes \alpha .g-v\otimes \alpha\in  \mathrm{Im}(\sigma)$. Whence $v\otimes \alpha.g-v\otimes \alpha\in V\otimes_L\omega$. In particular, for all $v\in V$, we must have $v\otimes (g-1)=v\otimes 1.g-v\otimes 1\in V\otimes \omega$. However $g-1$ belongs to $\F L$, so that $v\otimes (g-1)=v\otimes (g-1)1=v(g-1)\otimes 1\in V\otimes_L 1$. Thus $v(g-1)\otimes 1\in  V\otimes_L 1\cap V\otimes_L \omega=0$. But $v(g-1)\otimes 1$ can be uniquely written in the form 
$v_1\otimes 1+\sum_{x\in T}v_x\otimes (x-1)$ hence $v(g-1)=0$. Since  $v$ is a generic element in  the faithful $N$-module $V$, the element  $g$ must be $1$, as claimed.

We endow $M$ with a structure of $K$-module. The group $H$ acts by conjugation on $\F N$ and this action normalizes $\omega$ and all its powers, so that an $H$-action is induced on $\F N/\omega^k$. Given $v\otimes \alpha\in M$ and $g=hx\in K$, with $h\in H, \, x\in N$, define $v\otimes \alpha.g=vh\otimes \alpha^hx$. If $g=hx=h_1x_1$ we have
$vh_1\otimes \alpha^{h_1}x_1=v(hx{x_1}^{-1})\otimes \alpha^{h_1}x_1=
vh(x{x_1}^{-1})\otimes \alpha^{h_1}x_1$. Since $x{x_1}^{-1}=h^{-1}h_1\in N\cap H=L$, this element can be pulled through the tensor, so that 
$$vh(x{x_1}^{-1})\otimes \alpha^{h_1}x_1=
vh\otimes (x{x_1}^{-1})\alpha^{h_1}x_1= vh\otimes (h^{-1}h_1)\alpha^{h_1}x_1.$$
Now we easily get 
$$vh\otimes (h^{-1}h_1)\alpha^{h_1}x_1=vh\otimes (h^{-1})\alpha h (h^{-1}{h_1})x_1=
vh\otimes \alpha^ h x$$
thus showing that  the element $v\otimes \alpha.g$ is well defined. Now let $g=hx, \, g_1=h_1x_1$ be in $K$, and choose $v\otimes \alpha\in M$.  
We have $(v\otimes \alpha .g).g_1=(vh\otimes \alpha^hx)g_1=
v(hh_1)\otimes (\alpha^hx)^{h_1}x_1=v(hh_1)\otimes \alpha^{hh_1}x^{h_1}x_1$. On the other hand  
$v\otimes\alpha.(gg_1)=v\otimes\alpha.(hh_1^{x^{-1}}xx_1)=
v(hh_1^{x^{-1}})\otimes {\alpha}^{hh_1^{x^{-1}}}xx_1$. The element $z=[h_1, x^{-1}]$ belongs to $L$ and can be therefore pulled through the tensor symbol. Hence
$$v(hh_1^{x^{-1}})\otimes {\alpha}^{hh_1^{x^{-1}}}xx_1=v(hh_1z)\otimes {\alpha}^{hh_1z}xx_1
=v(hh_1)\otimes zz^{-1}{\alpha}^{hh_1}zxx_1.$$
Since $z=[h_1, x^{-1}]$ we have $zxx_1=x^{h_1}x_1$ thus showing that 
$$v\otimes\alpha.(gg_1)=v(hh_1)\otimes \alpha^{hh_1}x^{h_1}x_1=(v\otimes\alpha.g).g_1$$
for all $g, g_1\in K$ and $v\otimes \alpha$ in $M$. Since the elements $v\otimes \alpha$ generate $M$, we have an action of $K$ on $M$. We remark that the restriction of this action to $L$, gives the original structure of $L$-module on  $M$. Another important observation is the fact that $M$ is faithful as an $H$-module. To prove this fact we remind that we have already shown that the elments $v\otimes (1+\omega^k)$ are not $0$ in $M$, unless $v=0$. Let $g\in H$ be an element centralizing $M$. In particular $v\otimes (1+\omega^k)=(v\otimes (1+\omega^k)).g$ for all $v\in V$. Thus $v(g-1)\otimes (1+\omega^k)=0$ forcing $v(g-1)=0$ for all $v\in V$. Since $V$ is faithful this proves that $g=1$.
The action of $H$ on $M$ needs to be considered in more depth. Let $v\in V, \, a\in N$ and $h\in H$ and let us  calculate $[v\otimes (a+\omega^k), h]$. Using the identity $a^h=[h,a^{-1}]a$ and the fact that $[H,N]\leq L$, we get 
$$
v\otimes (a+\omega^k).h=
vh\otimes (a^h+\omega^k)=vh\otimes ([h,a^{-1}]a+\omega^k)=
vh^{a^{-1}}\otimes (a+\omega^k)
$$
and from this 
$$
[v\otimes (a+\omega^k),h]=[v,h^{a^{-1}}]\otimes (a+\omega^k)
$$
It is now easy to see that, given any set $h_1, \dots , h_s\in H$, we have 
$$
[v\otimes (a+\omega^k), h_1, \dots , h_s]=[v, {h_1}^{a^{-1}}, \dots , {h_s}^{a^{-1}}]\otimes (a+\omega^k)
$$
For any given $h\in H$ there exists $n\in \mathbb N$, such that $[V,_n h^G]=0$. Thus, for every $a\in N$, we have $[V\otimes_L(a+\omega^k),_n h^G]=0$. Since $M$ is generated, as a vector space, by  the set $\{v\otimes (a+\omega^k)\mid v\in V \, a\in N\}$,   the module  $M$ affords an $\mathcal F$-representation  for $H$. The preceeding argument shows that, if $h\in H$ and $[V,_n h^H]=0$, then $[M, _n h^G]=0$. 

 For every $t\in G$ define the $K$-module $M_t$ as follows. Take $M_t=M$ as abelian group and, for every $g\in K$ and $m\in M_t$, define $m.g= m.g^t$.  We let $K$ act on $\overline{M}=\mathrm{Dr}\{M_g\mid g\in G\}$ in the natural way. The dimension of $\overline{M}$ as a vector space over $\F$ is still bounded by $\left|G\right|$.  The kernel $P$  of the action of $K$ on $\overline{M}$ is normal in $G$ since $P=\mathrm{core}_G(C_K(M))$. Two cases should be considered.

\

{\bf Case 1} The group $P$ is non-trivial.

\

In this situation we consider the group $T=HP$. Since $H\cap P\leq C_H(M)=1$, the group $T$ is the direct product of $H$ and $P$. The group $P$ is nilpotent so that its center  is non-trivial. The group $H Z(P)$ is still normal in $G$ and strictly bigger than $H$, so we can find $A\leq Z(P)$ such that $A\trianglelefteq G$ and $A$ is torsionfree or elementary abelian of exponent $p$, according to $\F$ being $\mathbb Q$ or the field with $p$ elements. Given any $\F$-vector space  $U$ of dimension $\left| G\right|$, the group $A$ can be embedded into the $\F$-vector space $\mathrm{Hom}_\F(U,U)$, via some  homomorphism $\eta$. For each $a\in P$ define the $\F$-endomorphism $(a)\phi$ of $R=U\oplus U$ by setting $(u,w)(a)\phi=(u, w+u(a)\eta)$. It is clear that the map $\phi:P\longrightarrow \gl{R}{\F}$ is a faithful representation and $[r,_2 (a)\phi]=0$ for all $a\in P$ and $r\in R$. Since $P$ is normal in $G$, we have $[R, _2(a^G)\phi]=0$ for all $a\in P$, so that  $\phi$ is an $\
 mathcal F$-representation. The representation $\tau\otimes \phi$ on the tensor product $V\otimes_\F R$ is  a faithful $\mathcal F$-representation and, being $V\otimes_\F R$ of dimension at most $\left| G\right|$, any embedding of $V\otimes_\F R$ into $V$ gives rise to an $\mathcal F$-representation $\lambda$ of $T$ on $V$. As we have pointed out above, for each $h\in H$ we have $[V,_n (h^G)\tau]=0\, \Longrightarrow \, [V, (h^G)\lambda]=0$, proving   that $(T,\lambda)$ belongs to $\mathcal N$ and  $(H, \tau)\preccurlyeq (T, \lambda)$. On the other hand $H$ is properly contained in $T$, whence $(T, \lambda)\preccurlyeq (H, \tau)$ can not hold. This is a contradiction because $(H, \tau)$ was a nearly maximal element of $\mathcal N$

\bigskip

{\bf Case 2} The group $P$ is trivial.

\

The module $\overline{M}$ is then faithful for $K$ and can be embedded into $V$. Given such an embedding we get a faithful  $\mathcal F$-representation $\lambda$ of $K$ on $V$. As pointed out above $(K, \lambda)$ is in $\mathcal N$, $(H, \tau)\preccurlyeq (K, \lambda)$ but  $(K, \lambda)\not \preccurlyeq (H, \tau)$. This contradiction shows that  $H=G$ and the theorem is proved. \hfill $\Box$

\

We can now prove Theorem \ref{fitting as hp}

\pf First of all write $X$ as the union of an ascending chain of finite subsets $X_i$ for $ i\in \omega$, and set $N_i={X_i}^G$. By theorem \ref{rep} the group $G$ admits a faithful $\mathcal F$-representation over an $\F$-vector space $V$. We construct a $G$-series in $V$ inductively.
Start by setting $\s_0=\{ [V,_i N_0]\mid i\in \mathbb N\}\cup \{V\}$. The set $\s_0$ is a finite $G$-series in $V$, stabilized by $N_0$.
Assume the finite series $\s_{k}$ has been defined, in such a way that its elements are $G$-spaces and $N_k$ stabilizes $\s_{k}$.
We can now consider the action of $G$ on a factor $M$ of $\s_{k}$. The representation afforded by such factor is 
 still an $\mathcal F$-representations, so that the set $\{[M,_i N_{k+1}]\mid i\in \mathbb N\}$ is finite. The series $\s_{k+1}$ is obtained by adding to $\s_k$ all the preimages of the spaces $[M,_i N_{k+1}]$, for each factor $M$ of $\s_k$. The series $\s_{k+1}$ is a finite $G$-series and it is stabilized by $N_{k+1}$.
 The series $\s=\bigcup_{k\in \omega}\s_k$ is a series stabilized by $G$. Thus $G$ can be embedded into $S(\s)$. Since every element of $G$ lies in one of the $N_k$, each element stabilizes a finite subseries of $\s$. Therefore, given $g\in G$, the subgroup $g^{S(\s)}$ stabilizes a finite subseries too, so that it turns out to be nilpotent. Hence $G$ is contained in the Hirsch-Plotkin radical of $S(\s)$.\hfill $\Box$

It is not difficoult to find, for any given infinite cardinal  $\kappa$,  examples of Fitting groups $G$ of cardinality $\kappa$, of the form $G=X^G$ for some countable subset $X\subseteq G$. 

\bigskip

{\bf Example}  Let $\mathbb F$ be any field of cardinality $\kappa\geq \aleph_0$, and let 
$G=M(\mathbb Q,\mathbb F)$ be the McLain group with order type $\mathbb Q$. As usual, if 
$\mathcal B=\{v_q\mid q\in\mathbb Q\}$ is a basis for the natural module $V$ of $G$,   
we define the endomorphis $e_{rs}$ of $V$ by $v_qe_{rs}=\delta_{qr}v_s$. The group $G$ has cardinality $\kappa$ and it is an easy matter 
to show that, if $X=\{1+e_{rs}\mid r,s\in \mathbb Q, \,\, r<s\}$, then $X^G=G$.

\end{document}